\documentclass[11pt]{article}
\usepackage{amssymb}
\usepackage{amssymb,epsfig,latexsym}
\usepackage{graphicx,amsmath,amssymb,latexsym,epsfig}

\newtheorem{theorem}{Theorem}

\newtheorem{lemma}{Lemma}[theorem]
\newtheorem{proposition}{Proposition}[theorem]

\newtheorem{corollary}{Corollary}[theorem]
\newtheorem{remark}{Remark}
\newtheorem{example}{Example}

\def\EE{{\mathbb E}}
\def\PP{{\mathbb P}}

\def\tr{{\rm tr}}

\def\deg{{\mathrm{deg}}}

\def\ind{{\mathbf{1}}}
\def\bp{\noindent {\em Proof.} $\quad$}
\def\ep{\hfill $\Box$}
\def\rank{{\mathrm{rank}}}

\def\R{{\mathbb R}}
\def\C{{\mathbb C}}
\def\N{{\mathbb N}}
\def\cG{{\mathcal G}}

\def\cF{{\mathcal F}}
\def\cH{{\mathcal H}}
\def\cP{{\mathcal P}}

\def\Tr{{\mathrm{Tr}}}
\def\oG{\overline{G}}
\def\oX{\overline{X}}

\def\orho{\overline{\rho}}
\def\ocG*{\overline{\cG^*}}
\def\la{\langle}
\def\ra{\rangle}

\title{Resolvent of Large Random Graphs}
\author{Charles Bordenave \footnote{Institut de Math\'ematiques - Universit\'e de Toulouse \& CNRS  - France. Email: bordenave@math.univ-toulouse.fr } \, and \,  Marc Lelarge\footnote{INRIA-ENS - France. Email: marc.lelarge@ens.fr} }

\begin{document}

\maketitle

\begin{abstract}
We analyze the convergence of the spectrum of large random graphs to the spectrum of a limit infinite graph. We apply these results to graphs converging locally to trees and derive a new formula for the Stieltjes transform of the spectral measure of such graphs. We illustrate our results on the uniform regular graphs, Erd\"os-R\'enyi graphs and graphs with a given degree sequence. We give examples of application for weighted graphs, bipartite graphs and the uniform spanning tree of $n$ vertices.

\noindent {\em MSC-class:} 05C80, 15A52 (primary), 47A10 (secondary).
\end{abstract}


\section{Introduction}

Since the seminal work of Wigner \cite{wigner}, the spectral theory of
large dimensional random matrix theory has become a very active field
of research, see e.g. the monographs by Mehta \cite{mehta}, Hiai  and
Petz \cite{hiaipetz}, or Bai and Silverstein
\cite{bai-silverstein-book},  and for a review of applications in
physics, see Guhr et al. \cite{reviewRMT}. It is worth noticing that the classical random matrix theory has left aside the dilute random matrices (i.e. when the number of
non-zero entries on each row does not grow with the size of the
matrix). In the physics literature,
the analysis of dilute random matrices has been initiated by
Rodgers and Bray \cite{rodgersbray}. In \cite{bm99}, Biroli and
Monasson use heuristic arguments to anaylze the spectrum of the
Laplacian of Erd\"os-R\'eyni random graphs and an explicit connection
with their local approximation as trees is made by Semerjian and Cugliandolo
in \cite{sc02}.
Also related is the recent cavity approach to the
spectral density of sparse symmetric random matrices by Rogers et
al. \cite{r08}. 
Rigorous mathematical treatments can be found in Bauer and Golinelli
\cite{bauer01} and Khorunzhy, Scherbina and Vengerovsky
\cite{khorunzhy04} for Erd\"os-R\'eyni random graphs. 
In parallel, since McKay \cite{mckay}, similar
questions have also appeared in graph theory and combinatorics, for a
review, refer to Mohar and Weiss \cite{mohwoe}. In this paper, we
present a unified  treatment of these issues, and prove under weak
conditions the convergence of the empirical spectral distribution of
adjacency and Laplacian matrices of large graphs.

Our main contribution is to connect this convergence to the local weak
convergence of the sequence of graphs. There is a growing interest in
the theory of convergence of graph sequences. The convergence of dense
graphs is now well understood thanks to the work of Lov\'asz and Szegedy
\cite{ls06} and a series of papers written with Borgs, Chayes S\'os
and Vesztergombi (see \cite{bc1,bc2} and references therein).
They introduced several natural metrics for graphs and showed that
they are equivalent. However these results are of no help in the case
studied here of diluted graphs, i.e. when the number of edges scales
as the number of vertices. Many new phenomena occur, and there are a
host of plausible metrics to consider \cite{br08}. Our first main
result (Theorem \ref{th:convmes}) shows that the local weak convergence implies the convergence
of the spectral measure.
Our second main result (Theorem \ref{th:adj}) characterizes in term of
Stieltjes transform the limit spectral measure of a large class of
random graphs ensemble.
The remainder of this paper is organized as follows: in the next
section, we give our main results. In Section
\ref{sec:convop}, we prove Theorem \ref{th:convmes}, in Section
\ref{sec:Th2} we prove Theorem \ref{th:adj}. Finally, in Section
\ref{sec:ext}, we extend and apply our results to related graphs.
\

\section{Main results}\label{sec:main}

\subsection{Convergence of the spectral measure of random graphs}\label{subsec:conv}
Let $G_n$ be a sequence of simple graphs with vertex set $[n]=\{1,\dots,
n\}$ and undirected edges set $E_n$.
We denote by $A =  A(G_n)$, the $n \times n$ adjacency
matrix of $G_n$, in which $A_{ij}=1$ if $(ij) \in E_n $ and
$A_{ij}=0$ otherwise. The Laplace matrix of $G_n$ is
$L(G_n)= D(G_n)-A(G_n)$, where $D=D(G_n)$ is the degree diagonal
matrix in which $D_{ii}=\deg(G_n,i) := \sum_{j \in [n]}
A_{ij}$ is the degree of $i$ in $G_n$ and $D_{ij}=0$
for all $i\neq j$. The main object of this paper is to study the
convergence of the empirical measures of the eigenvalues of
$A$ and $L$ respectively when the sequence of graphs
converges weakly as defined by Benjamini and Schramm \cite{bensch} and
Aldous and Steele \cite{aldste} (a precise definition is given below). 
Note that the spectra of $A(G_n)$ or $L(G_n)$
do not depend on the labeling of the graph $G_n$. If we label the
vertices of $G_n$ differently, then the resulting matrix is
unitarily equivalent to $A(G_n)$ and $L(G_n)$ and it is
well-known that the spectra are unitarily invariant. For ease of
notation, we define $$ \Delta_n   = A(G_n) - \alpha
D(G_n),$$ with $\alpha \in \{0,1\}$ so that $\Delta_n =
A(G_n)$ if $\alpha = 0$ and  $\Delta_n = - L(G_n)$ if
$\alpha = 1$. The empirical spectral measure of $\Delta_n$ is denoted
by $$\mu_{\Delta_n}   = n^{-1} \sum_{i=1}^n \delta_{\lambda_{i}(\Delta_n)},$$
where $(\lambda_{i}(\Delta_n))_{1 \leq i \leq n}$ are the eigenvalues
of $\Delta_n$. 
We endow the set of measures on $\R$ with the
usual weak convergence topology. This convergence is metrizable
with the L\'evy distance $L(\mu,\nu) = \inf \{ h \geq 0 : \forall x
\in \R , \mu((-\infty,x-h])-h \leq \nu ((-\infty,x]) \leq
\mu((-\infty,x-h]) + h \}$.

We now define the local weak convergence introduced by Benjamini and
Schramm \cite{bensch} and Aldous and Steele \cite{aldste}.
For a graph $G$, we define the rooted graph $(G,o)$ as the connected
component of $G$ containing a distinguished vertex $o$ of $G$, called
the root. 
A homomorphism form a graph $F$ to another graph $G$ is an
edge-preserving map form the vertex set of $F$ to the vertex set of
$G$.
A bijective homomorphism is called an isomorphism.
A rooted isomorphism of rooted graphs is an isomorphism of the graphs
that takes the root of one to the root of the other. 
$[G,o]$ will denote the class of rooted graphs that are
rooted-isomorphic to $(G,o)$.
Let $\cG^*$ denote the set of rooted isomorphism classes of rooted
connected locally finite graphs. Define a metric on $\cG*$
by letting the distance between $(G_1,o_1)$ and $(G_2,o_2)$ be
$1/(\alpha+1)$, where $\alpha$ is the supremum of those $r\in \N$ such
that there is some rooted isomorphism of the balls of radius $r$ (for
the graph-distance) around the roots of $G_i$.
$\cG^*$ is a separable and complete metric space \cite{aldlyo}.
For probability measures $ \rho, \rho_n$ on
$\cG^*$, we write $\rho_n\Rightarrow  \rho$
when $\rho_n$ converges weakly with respect to this
metric.

Following \cite{aldlyo}, for a finite graph $G$, let $U(G)$ denote the
distribution on $\cG^*$ obtained by choosing a uniform random vertex of
$G$ as root.
We also define $U_2(G)$ as the distribution on $\cG^* \times \cG^*$ of
the pair of rooted graphs $((G,o_1),(G,o_2))$ where $(o_1,o_2)$ is a
uniform random pair of vertices of $G$.
If $(G_n), n \in \N,$ is a sequence of {\em deterministic}
graphs with vertex set $[n]$
and $\rho$ is a probability measure on $\cG^*$, we say
the random weak limit of $G_n$ is $\rho$ if $U(G_n)
\Rightarrow \rho$. 
If $(G_n), n \in \N,$ is a sequence of {\em random} graphs with vertex
set $[n]$, we
denote by $\EE[\cdot]=\EE_n[\cdot]$ the expectation with respect to the randomness
of the graph $G_n$. The measure $\EE[U(G_n)]$ is defined as
$\EE[U(G_n)](B)=\EE[U(G_n)(B)]$ for any measurable event $B$ on $\cG^*$.
Following Aldous and Steele \cite{aldste}, we will say that
the random weak limit of $G_n$ is $\rho$ if $ \EE[U(G_n)]\Rightarrow
\rho$.
Note that the second definition generalizes the first one (take
$\EE_n= \delta_{G_n}$).
In all cases, we denote by $(G,o)$ a random rooted graph whose
distribution of its equivalence class in $\cG^*$ is $\rho$.
Let $\deg(G_n,o)$ be the degree of the root under
$U(G_n)$ and $\deg(G,o)$ be the degree of the root under $\rho$. They are random variables on $\N$ such that if the random weak limit
of $G_n$ is $\rho$ then $\lim_{n\to\infty}\EE[\deg(G_n,o)\leq
k]=\rho(\{\deg(G,o)\leq k\})$.

We will make the following assumption for the whole paper:
\begin{itemize}
\item[A.] The sequence of random variables $(\deg (G_n,o)), n \in \N,$ is uniformly integrable.
\end{itemize}

Assumption (A) ensures that if the random weak limit of $G_n$ is $\rho$,
then the average degree of the root converges, namely
$
\lim_{n\to\infty}\EE[\deg(G_n,o)] = \rho(\deg(G,o)).
$

To prove our first main result, we will consider two assumptions,
one, denoted by (D), for a given sequence of finite graphs  and
another,  denoted by (R), for a sequence of random finite graphs.
 \begin{itemize}
\item[D.] As $n$ goes to infinity, the random weak limit of $G_n$ is $\rho$.
\item[R.] As $n$ goes to infinity, $U_2(G_n)\Rightarrow \rho \otimes \rho$.
\end{itemize}

Of course, Assumption (R) implies (D). We are now ready to state our first
main theorem:

\begin{theorem}\label{th:convmes}
\begin{itemize}
\item[(i)] Let $G_n = ( [n] , E_n )$ be a sequence of graphs
satisfying assumptions (D-A), then there exists a probability
measure $\mu  $ on $\R$ such that $\lim_{n \to \infty} \mu_{\Delta_n}
= \mu  $. 
\item[(ii)] Let $G_n = ( [n] , E_n )$ be a sequence of
random graphs satisfying assumptions (R-A), then there exists a
probability measure $\mu  $ on $\R$ such that, $\lim_{n \to
\infty} \EE L (\mu_{\Delta_n}  , \mu  ) = 0$.
\end{itemize}
\end{theorem}

In (ii), note that the stated convergence implies the weak
convergence  of the law of $\mu_{\Delta_n}$ to $\delta_{\mu  }$.
Theorem \ref{th:convmes} appeared under different settings, when the
sequence of maximal degrees of the graphs $G_n$ is bounded, see Colin
de Verdi\`ere \cite{colindeverdiere}, Serre \cite{serre} and Elek
\cite{elek}.

\subsection{Random graphs with trees as local weak limit}

\label{subsec:tree}

\begin{figure}[htb]
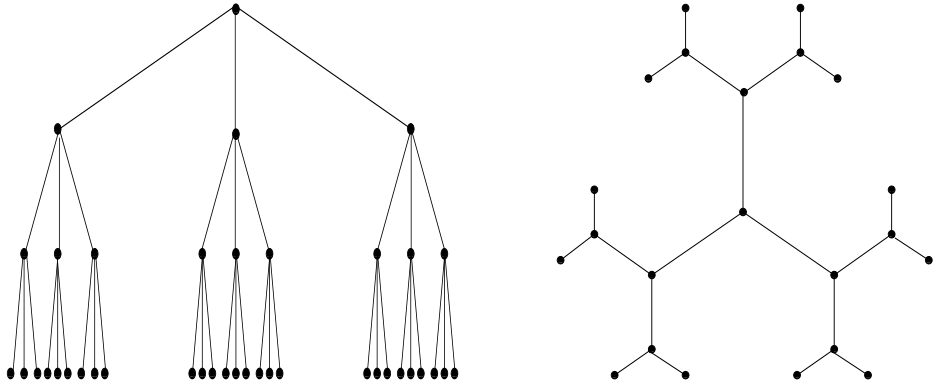

\begin{center}
 \includegraphics[angle=0,width = 6cm, height=5cm]{3ary.eps} \hspace{30pt}
 \includegraphics[angle=0,width = 5cm, height=5cm]{3reg.eps}
\caption{Left: representation of a $3$-ary tree. Right:
representation of a $3$-regular tree.} \label{fig1}
\end{center}\end{figure}

We now consider a sequence of random graphs $G_n, n \in \mathbb
N$ which converges as $n$ goes to infinity to a possibly
infinite tree. In this case, we will be able to characterize the
probability measure $\mu$.
Here, we restrict our attention to particular
trees as limits but some cases outside the scope of this section are
also analyzed in Section \ref{sec:ext}.
A Galton-Watson Tree (GWT) with {\em offspring}
distribution $F$ is the random tree obtained by a standard
Galton-Watson branching process with offspring distribution $F$.
For example, the infinite $k$-ary tree is a GWT with offspring
distribution $\delta_k$,  see Figure \ref{fig1}. A GWT with {\em
degree} distribution $F_*$ is a rooted random tree obtained by a
Galton-Watson branching process  where the root has offspring
distribution $F_*$ and all other genitors have offspring
distribution $F$ where for all $k \geq 1$, $F(k-1) = k F_* (k) /
\sum_k k F_* (k)$ (we assume $\sum_k k F_* (k) < \infty$). For
example the infinite $k$-regular tree is a GWT with degree
distribution $\delta_k$, see Figure \ref{fig1}.  It is easy to
check that a GWT with degree distribution $F_*$ defines a
unimodular probability measure on $ \cG^*$ (for a definition and
properties of unimodular measures, refer to \cite{aldlyo}). Note
that if $F_*$ has a finite second moment then $F$ has a finite
first moment.

Let $n \in \N$ and $G_n = ( [n] , E_n )$, be a random graph on
the finite vertex set $[n]$ and edge set $E_n$. 
We assume that the following holds
\begin{itemize}
\item[RT.] As $n$ goes to infinity, $U_2(G_n)$ converges weakly to
$\rho \otimes \rho$, where $\rho \in \cG^*$ is the probability
measure of GWT with degree distribution $F_*$.
\end{itemize}

Note that we have $\lim_{n\to\infty} \EE[\deg(G_n,o)] = \sum_k kF_*(k)<\infty$
and $$\lim_{n\to\infty}\EE[\deg(G_n,v)|\:(o,v)\in E_n]=1+
\sum_kkF(k)=1+\sum_kk(k-1)F_*(k)/\sum_k kF_*(k).$$
Under assumption (A), this last quantity might be infinite. To
prove our next result, we need to strengthen  it into
\begin{itemize}
\item[A'.] Assumption (A) holds and $\sum_k k^2F_*(k)<\infty$.
\end{itemize}

We mention three important classes of graphs which converge
locally to a tree and which satisfy our assumptions.

\begin{example}{Uniform regular graph.}
\label{ex:uniform} \textnormal{ The uniform $k$-regular graph on
$n$ vertices satisfies these assumptions with the infinite
$k$-regular tree as local limit. It follows for example easily
from Bollob{\'a}s \cite{bollobas80}, see also the survey Wormald
\cite{wormald}. }
\end{example}

\begin{example}{Erd\"os-R\'enyi graph.}
\label{ex:erdosrenyi} \textnormal{ Similarly, consider the
Erd\"os-R\'enyi graph on $n$ vertices where there is an edge
between two vertices with probability $p/n$ independently of
everything else. This sequence of random graphs satisfies the
assumptions with limiting tree the GWT with degree distribution
$Poi(p)$.}
\end{example}

\begin{example}{Graphs with asymptotic given degree.}
\textnormal{ More generally, the usual random graph (called
  configuration model) with asymptotic degree distribution $F_*$
  satisfies this set of hypothesis provided that $\sum k(k-1) F_*(k) <
  \infty$ (e.g. see Chapter 3 in Durrett \cite{durrett2007} and Molloy and Reed
\cite{molloyreed}).}
\end{example}

In these three cases, it is easy to see that if $G_n$ denotes the
random graph on $[n]$ and $G$ is a GWT with degree distribution $F_*$,
then we have convergence of $G_n$ to $G$ in probability which implies
Assumption (RT).

Recall that $ \Delta_n  = A(G_n) - \alpha D(G_n),$ with
$\alpha \in \{0,1\}$.
We now introduce a standard tool of random matrix theory used to
describe the empirical spectral measure $\mu_{\Delta_n}$ (see Bai and
Silverstein \cite{bai-silverstein-book} for more details).
Let $R_n(z) = (\Delta_n -zI_n)^{-1}$ be the resolvent of $ \Delta_n$.
We denote $\C_+ = \{ z \in \C :
\Im z > 0 \}$. Let $\cH$ be the set of holomorphic functions $f$
from $\C_+$ to $\C_+$ such that $|f(z)|\leq \frac{1}{\Im z}$.
For all $i \in [n]$, the mapping $z \mapsto  R_n (z)_{ii}$ is in $\cH$
(see Section \ref{sec:rfn}).
We denote by $\cP(\cH)$ the space of
probability measures on $\cH$.
The Stieltjes transform of the empirical
spectral distribution is given by:
\begin{eqnarray}
\label{eq:sti}m_n(z)&=& \int_{\R} \frac{1}{x- z} d\mu_{\Delta_n}
(x)=\frac{1}{n} \tr R_n (z)=\frac{1}{n} \sum_{i=1}^n R_n(z)_{ii} 
\end{eqnarray}
where $z\in \C_+$.
Our main second result is the following.
\begin{theorem}
\label{th:adj} Under assumptions (RT-A'),
\begin{itemize}
\item[(i)] There exists a unique probability measure $Q   \in
\cP(\cH)$ such that for all $z\in  \mathbb C_+$,
\begin{equation}
\label{eq:RDE} Y (z)   \stackrel{d}{=}  - \left(z  + \alpha (N+1)
+ \sum_{i=1}^ {N} Y_i(z)  \right)^{-1},
\end{equation}
where $N$ has distribution $F$ and $Y$ and $Y_i$ are iid copies
independent of $N$ with law $Q  $. \item[(ii)] For all $z\in
\C_+$, $m_n   (z)$ converges as $n$ tends to infinity in $L^1$
to $ \EE X   (z)$, where for all $z \in \C_+$,
\begin{equation}
\label{eq:X} X   (z)\stackrel{d}{=}  - \left(z  + \alpha N_*+
\sum_{i=1}^ {N_*} Y_i(z)   \right)^{-1},
\end{equation}
where $N_*$ has distribution $F_*$ and $Y_i$ are iid copies with
law $Q  $, independent of $N_*$.
 \end{itemize}
\end{theorem}

Equation (\ref{eq:RDE}) is a Recursive Distributional Equation
(RDE). 
In random
matrix theory, the Stieltjes transform appears classically as a
fixed point of a mapping on $\cH$. For example, in the Wigner case \cite{wigner}
(i.e. the matrix $W_n = (A_{ij}/\sqrt n )_{1 \leq i, j \leq n}$
where $A_{ij} = A_{ji}$ are iid copies of $A$ with $var(A) = \sigma^2$),
the Stieltjes tranform $m(z)$ of the limiting spectral measure
satisfies for all $z \in \C_+$,
\begin{equation} \label{eq:fpwigner}
 m(z) = - \left( z + \sigma^2 m(z) \right)^{-1}. \end{equation}
It is then easy to show that the limiting spectral measure is
the semi circular law with radius $2\sigma$.

We explain now why the situation is different in our case.
Let $X_n(z)=R_n(z)_{oo}\in \cH$ be the
diagonal term of the resolvent matrix corresponding to the root of the
graph $G_n$.
From (\ref{eq:sti}), we see that $m_n(z)=U(G_n)(X_n(z))$. 
In a first step, we will prove that the random variable
$X_n(z)$ under $U(G_n)$ converges in distribution to $X(z)$ given by
(\ref{eq:X}).
Then we finish the proof of Theorem \ref{th:adj} with the following steps $\lim_n
m_n(z)=\lim_n\EE[U(G_n)(X_n(z))] = \EE[X(z)]$.
For the proof of the convergence of $X_n(z)$, we first consider the the case of uniformly bounded degrees in $G$, i.e. $\max_i\deg
(G_n,i)\leq \ell$.
In this case, Mohar proved in \cite{moh82} that it is possible to
define the resolvent $R$ of the possibly infinite graph $(G,o)$ with
law $\rho$. We then show that $X_n(z)$ converges weakly to
$R(z)_{oo}$ under $\rho$. Thanks to the tree structure of $(G,o)$, we
can characterize the law of $R(z)_{oo}$ with the RDE
(\ref{eq:RDE}). Recall that the offspring distribution of the root has
distribution $F_*$ whereas all other nodes have offspring distribution
$F$ which explains the formula (\ref{eq:X}) and (\ref{eq:RDE})
respectively.
A similar approach was used in \cite{bordenave-2009} for the spectrum
of large random reversible Markov chains with heavy-tailed weights,
where a more intricate tree structure appears.

Note that for all $z \in \C_+$, $m_n(z)$ and $X (z)$ are
bounded by $\Im (z)^{-1}$, hence the convergence in Theorem
\ref{th:adj}(ii) of $m_n(z)$ to $\EE X(z)$ holds in
$L^p$ for all $p\geq 1$. 
Under the restrictive assumption $\max_i\deg (G_n,i)\leq \ell$, we are
able to prove that on $\cH$, $X_n$ converges weakly to $X$ but we do
not know if this convergence holds in general.
We also need Assumption (A') to prove the uniqueness of the solution
in (\ref{eq:RDE}) even though we know from Theorem \ref{th:convmes}
that the empirical spectral distribution converges under the weaker
Assumption (A).

We end this section with two examples that appeared in the literature.

\noindent{ \textsc{Example} \ref{ex:uniform}} $\;$ If $G_n$ is the
uniform $k$-regular graph on $[n]$, with $k\geq 2$, then $G_n$
converges to the GWT with degree distribution $\delta_{k}$. We
consider the case $\alpha = 0$, looking for deterministic
solutions of $Y$, we find: $ Y (z) = - ( z + (k-1) Y(z))^{-1}$,
hence, in view of (\ref{eq:fpwigner}), $Y$ is simply the Stieltjes
transform of the semi-circular law with radius $2 \sqrt{k-1}$. For
$X(z)$ we obtain,
\begin{eqnarray}
\label{RDEreg}X (z) =  - ( z + k Y(z) )^ {-1} = - \frac{ 2 ( k-1)
}{
  (k-2) z + k \sqrt { z^2 - 4 (k-1) } }.
\end{eqnarray}
In particular $\Im X(z)  = \Im (z + k Y(z) ) / |z + k Y(z) |^2$.
Using the formula $\mu[a,b] = \lim_{v \to 0^+} \frac{1}{\pi}
\int_a ^ b \Im X(x + i v) dx$, valid for all continuity points
$a<b$ of $\mu$, we deduce easily  that $\mu_{A_n}$ converges
weakly to the probability measure $\mu(dx) = f(x) dx$ which has a
density $f$ on $[-2 \sqrt{k-1}, 2 \sqrt{k-1}]$ given by
$$
f(x) =  \frac{k }{2\pi} \frac{ \sqrt { 4 (k-1) - x^2 } }{ k^2 -
x^2 },
$$
and $f(x) = 0$ if $x \notin [-2 \sqrt{k-1}, 2 \sqrt{k-1}]$. This
formula for the density of the spectral measure is due to McKay
\cite{mckay} and Kesten \cite{kesten59} in the context of simple
random walks on groups. To the best of our knowledge, the proof of
Theorem \ref{th:adj} is the first proof using the resolvent method
of McKay's Theorem. It is interesting to notice that this measure
and the semi-circle distribution are simply related by their
Stieltjes transform see (\ref{RDEreg}).

\noindent{{\sc Example \ref{ex:erdosrenyi}}} $\;$ If $G_n$ is a
Erd\"os-R\'enyi graph on $[n]$, with parameter $p/n$ then $G_n$
converges to the GWT with degree distribution $Poi(p)$. In this
case, $F$ and $F_*$ have the same distribution,  thus for
$\alpha = 0$, $X \stackrel{d}=Y$ has law $Q$. Theorem \ref{th:adj} improves
a result of Khorunzhy, Shcherbina and Vengerovsky
\cite{khorunzhy04}, Theorems 3 and 4, Equations (2.17), (2.24).
Indeed, for all $z \in \C_+$, we use the formula $
e^{i  u  w} = 1- \sqrt u \int_{0} ^ {\infty} \frac{ J_1
( 2 \sqrt {u t} )}{\sqrt t} e^ {  - i t w^{-1} }   dt $ valid
for all $u \geq 0 $ and $w \in \C_+$, and where $J_1 ( t ) =
\frac{t}{2} \sum_{k = 0 } ^ {\infty} \frac{( - t ^2 / 4) ^ k }{k !
  (k+1)! }$ is the Bessel function of the first kind. Then with $f(u,
z) = \EE e^{Êu Y(z)} $ and $\varphi (z) = \EE z^N$ we obtain easily
from (\ref{eq:RDE}) that for all $u \geq 0$, $z \in \C_+$,
  $$
  f(u,z) =  1 - \sqrt u \int_{0} ^ {\infty} \frac{ J_1
( 2 \sqrt {u t} )}{\sqrt t} e^ {  i t ( z + \alpha) } \varphi( e^{ i \alpha t } f( t, z) )    dt.
  $$
  If $N$ is a Poison random variable with parameter $p$, then $\varphi (z) = \exp ( p ( z - 1) )$, and we obtain the results in \cite{khorunzhy04}.

\section{Proof of Theorems \ref{th:convmes}}\label{sec:convop}

\subsection{Random finite networks}\label{sec:rfn}

It will be convenient to work with marked graphs that we call networks.
First note that the space $\cH$ of holomorphic functions $f$
from $\C_+$ to $\C_+$ such that $|f(z)|\leq \frac{1}{\Im z}$ equipped
with the topology induced by the uniform convergence on compact sets
is a complete separable metrizable compact space (see Chapter 7 in
\cite{con}).
We now define a network as a graph $G=(V,E)$ together with a complete separable
metric space, in our case $\cH$, and a map from $V$ to $\cH$.
We use the following notation: $G$ will denote a graph and $\oG$ a network
with underlying graph $G$. A rooted network $(\oG,o)$ is the connected
component of a network $\oG$ of a distinguished vertex $o$ of $G$,
called the root.
$[\oG,o]$ will denote the class of rooted networks that are
rooted-isomorphic to $(\oG,o)$.
Let $\ocG*$ denote the set of rooted isomorphism classes of rooted
connected locally finite networks. As in \cite{aldlyo}, define a
metric on $\ocG*$ by letting the distance between $(\oG_1,o_1)$ and
$(\oG_2,o_2)$ be $1/(\alpha+1)$, where $\alpha$ is the supremum of
those $r\in \N$ such that there is some rooted isomorphism of the
balls of (graph-distance) radius $r$ around the roots of $G_i$ such
that each pair of corresponding marks has distance less than $1/r$.
Note that the metric defined on $\cG^*$ in Section \ref{subsec:conv}
corresponds to the case of a constant mark attached to each vertex.
It is now easy to extend the local weak convergence of graphs to
networks.
To fix notations, for a finite network $\oG$, let $U(\oG)$ denote the
distribution on $\ocG*$ obtained by choosing a uniform random vertex
of $G$ as root.
If $(\oG_n), n \in \N,$ is a sequence of (possibly random) networks
with vertex set $[n]$, we denote by $\PP_n[\cdot]$ the expectation
with respect to the randomness of the graph $G_n$. 
We will say that the random weak limit of $\oG_n$ is $\orho$ if $
\PP_n[U(\oG_n)]\Rightarrow \orho$.

In this paper, we consider the finite networks
$\{G_n, (R_n(z)_{ii})_{i\in [n]}\}$, where we attach the mark
$R_n(z)_{ii}$ to vertex $i$. We need to check that the map $z\mapsto
R_n(z)_{ii}$ belongs to $\cH$.
First note that by standard linear algebra (see Lemma \ref{le:schur} in Appendix), we
have
\begin{eqnarray*}
R_n(z)_{ii} =
\frac{1}{(\Delta_n)_{ii}-z-\beta_i^t(\Delta_{n,i}-zI_{n-1})^{-1}\beta_i},
\end{eqnarray*}
where $\beta_i$ is the $((n-1)\times 1)$ $i$th column vector of
$\Delta_n$ with the $i$th element removed and $\Delta_{n,i}$
is the matrix obtained form $\Delta_n$ with the $i$th row and
column deleted (corresponding to the graph with vertex $i$
deleted). An easy induction on $n$ shows that $R_n(z)_{ii}\in
\cH$ for all $i\in [n]$.

\subsection{Linear operators associated with a graph of bounded degree}\label{sec:graphop}

We first recall some standard results that can be found in Mohar
and Woess \cite{mohwoe}. Let $(G,o)$ be the connected component of a
locally finite graph $G$ containing the vertex $o$ of $G$.
There is no loss of generality in assuming that the vertex set of $G$
is $\N$.
Indeed if $G$ is finite, we can extend $G$ by adding isolated vertices.
We assume that $\deg(G)=\sup\{\deg(G,u),
u\in \N\}<\infty$. Let $A(G)$ be the adjacency matrix of $G$. We
define the matrix $\Delta=\Delta(G) = A(G)-\alpha D(G)$, where $D(G)$
is the degree diagonal matrix and $\alpha\in \{0,1\}$. Let
$e_k=\{\delta_{ik}: i\in \N\}$ be the specified complete
orthonormal system of $L^2(\N)$. Then $\Delta$ can be interpreted
as a linear operator over $L^2(\N)$, which is defined on the basis
vector $e_k$ as follows:
\begin{eqnarray*}
\la \Delta e_k, e_j\ra = \Delta_{kj}.
\end{eqnarray*}
Since $G$ is locally finite, $\Delta e_k$ is an element of
$L^2(\N)$ and $\Delta$ can be extended by linearity to a dense
subspace of $L^2(\N)$, which is spanned by the basis vectors
$\{e_k,k\in \N\}$. Denote this dense subspace $H_0$ and the
corresponding operator $\Delta_0$. The operator $\Delta_0$ is
symmetric on $H_0$ and thus closable (Section VIII.2 in
\cite{reesim}). We will denote the closure of $\Delta_0$ by the
same symbol $\Delta$ as the matrix. The operator $\Delta$ is by
definition a closed symmetric transformation: the coordinates of
$y=\Delta x$ are
\begin{eqnarray*}
y_i = \sum_{j}\Delta_{ij}x_j,\quad i\in \N,
\end{eqnarray*}
whenever these series converge.

The following lemma is proved in \cite{moh82} for the case
$\alpha=0$ and the case $\alpha=1$ follows by the same argument.
\begin{lemma}
Assume that $\deg(G)=\sup\{\deg(G,u), u\in V\}<\infty$, then
$\Delta$ is self-adjoint.
\end{lemma}

\subsection{Proof of Theorem  \ref{th:convmes} (i) - bounded degree}

\label{subsec:th(i)bounded}

In this paragraph, we assume that there exists $\ell\in \N$ such that
\begin{eqnarray}
\label{eq:bdeg}\rho(\deg(G,o)\leq \ell)=1.
\end{eqnarray}

Since $\cG^*$ is a complete separable metric space, by the Skorokhod
Representation Theorem (Theorem \ref{th:skorokhod} in Appendix), we
can assume that $U(G_n)$ and $(G,o)$ are defined on a common
probability space such that $U(G_n)$ converges almost surely to
$(G,o)$ in $\cG^*$.
As explained in previous section, we define the operators $\Delta =
\Delta(G)$ and $\Delta_n = \Delta(G_n)$ on the Hilbert space $L^2
(\N)$. 
The convergence of $U(G_n)$ to $(G,o)$ implies a convergence of the
associated operators up to a re-indexing of $\N$ which preserves the
root.
To be more precise, let $(G,o)[r]$ denote the finite graph induced by
the vertices at distance at most $r$ from the root $o$ in $G$.  
Then by definition there exists a.s. a sequence $r_n$ tending to
infinity and a random rooted isomorphism $\sigma_n$ from $(G,o)[r_n]$
to $U(G_n)[r_n]$. 
We extend arbitrarily this isomorphism $\sigma_n$ to all $\N$. 
For the basis vector $e_k \in L^2 (\N)$, we set $\sigma_n ( e_k )  = e
_ {\sigma_n (k) }$. For simplicity, we assume that $e_0$ corresponds
to the root of the graph (so that $\sigma_n(0)=0$) and we denote
$e_0=o$ to be consistent with the notation used for graphs.
By extension, we define $\sigma_n(  \phi ) $ for
each $\phi\in H_0$ the subspace of $L^2(\N)$, which is spanned by the
basis vectors $\{e_k,k\in \N\}$. 
Then the convergence of $U(G_n)$ to $(G,o)$ implies that for all $\phi
\in H_0$,
\begin{eqnarray}
\label{conv} \sigma_n ^{-1} \Delta_n \sigma_n \phi\to \Delta\phi \mbox{ a.s.}
\end{eqnarray}
By Theorem VIII.25(a) in
\cite{reesim}, the convergence (\ref{conv}) and the fact that
$\Delta$ is a self-adjoint operator (due to (\ref{eq:bdeg})) imply the convergence of
$\sigma_n^{-1} \Delta_n \sigma_n \to\Delta$ in the strong resolvent sense:
\begin{eqnarray*}
\sigma_n ^{- 1} R_n(z) \sigma_n  x - R(z)x \to 0,\mbox{ for any $x\in
  L^2(\N)$, and for all $z\in \C_+$.}
\end{eqnarray*}
This last statement shows that the sequence of networks $U(\oG_n)$
converges a.s. to $(\oG,o)$ in $\ocG*$.
In particular, we have
\begin{eqnarray}
\label{eq:convbounded}m_n(z) =  \frac{1}{n} \sum_{i=1}^n \la R_n(z)e_i,e_i\ra = U(G_n)(\la R_n(z)o,o\ra)\to \int \la R(z)o,o\ra d\rho[G,o],
\end{eqnarray}
by dominated convergence since $|\la R_n(z)e_i,e_i\ra|\leq (\Im z)^{-1}$.  
\begin{remark}
A trace operator $\Tr$ was defined in \cite{aldlyo}. With their notation, we have
$$
 \Tr(R(z))=  \int \la R(z)o,o \ra d\rho[G,o].
$$
\end{remark}

\subsection{Proof of Theorem  \ref{th:convmes} (i) - general case} \label{sec:th(i)}

Let $\ell \in \N$, we define the graph $G_{n,\ell}$ on $[n]$
obtained from $G_n$ by removing all edges adjacent to a vertex
$i$, if $\deg(G_n,i) > \ell$.
Therefore the matrix $\Delta(G_{n,\ell})$ denoted by
$\Delta_{n,\ell}$ is equal to, for $i \neq j$
$$
 (\Delta_{n,\ell})_{ij} = \left\{\begin{array}{ll}
A(G_n)_{ij} & \hbox{ if }\max \{ \deg(G_n,i),\deg(G_n,j) \} \leq \ell  \\
0 & \hbox{ otherwise}
\end{array}\right.
$$
and  $(\Delta_{n,\ell})_{ii} = -\alpha \sum_{j \neq
i}  (\Delta_{n,\ell})_{ij}$. The empirical measure of the
eigenvalues of $  \Delta_{n,\ell}$ is denoted by
$\mu_{\Delta_{n,\ell}}$. By the Rank Difference Inequality (Lemma \ref{le:rankineq} in Appendix),
$$
L( \mu_{\Delta_n}, \mu_{\Delta_{n,\ell}}  ) \leq \frac{1}{n} \rank (
\Delta_{n} -  \Delta_{n,\ell} )  ,
$$
where $L$ denotes the L\'evy distance. The rank of $\Delta_n -
\Delta_{n,\ell}$ is upper bounded by the number of rows different from $0$, i.e. by the number of vertices with degree at least $\ell +1$ or such that there exist a neighboring vertex with degree at least $\ell +1$. By definition, each vertex with degree $d$ is connected to $d$ other vertices. It follows
$$
\rank (
\Delta_{n} -  \Delta_{n,\ell} )   \leq \sum_{i = 1} ^n ( \deg(G_n,i) +1)  \ind ( \deg(G_n,i)  > \ell),
$$
and therefore:
$$
L (\mu_{\Delta_n}, \mu_{\Delta_{n,\ell}}) \leq  \int ( \deg(G,o) +1) \ind ( \deg(G,o)  > \ell) dU(G_n) [G,o]=:p_{n,\ell}.
$$
By assumptions (D-A), uniformly in $\ell$,
\begin{eqnarray*}
p_{n,\ell}\to \int ( \deg(G,o) +1 ) \ind (
\deg(G,o)  > \ell)d\rho[G,o],
\end{eqnarray*}
where the right-hand term tends to $0$ as $\ell$ tends to infinity. Fix $\epsilon>0$, for $\ell$ sufficiently large and $n\geq N(\ell)$, we
have
$$
L( \mu_{\Delta_n} , \mu_{\Delta_{n,\ell}})  \leq \epsilon.$$
Hence we get, for $n\geq N(\ell)$, $q \in \N$,
\begin{eqnarray*}
L( \mu_{\Delta_{n+q}} ,  \mu_{\Delta_n} )  \leq  2 \epsilon +  L( \mu_{\Delta_{n+q,\ell}} , \mu_{\Delta_{n,\ell}}).
\end{eqnarray*}
By (\ref{eq:convbounded}), if $m_{n,\ell} (z) $ denotes the Stieltjes transform of
$\mu_{\Delta_{n,\ell}}$, we have for all $z \in \C_+$, $\lim_{n\to\infty}m_{n,\ell}
(z)=m^\ell(z)$ for some $m^\ell\in \cH$. Hence for the weak convergence $\lim_{n\to\infty} \mu_{\Delta_{n,\ell}}=\mu^\ell$ for the probability measure $\mu^\ell$ whose Stieltjes transform is $m^\ell$. It follows that the sequence $\mu_{\Delta_n}$
is Cauchy and the proof of Theorem \ref{th:convmes}(i) is
complete (recall that the set probability measures on $\R$ with the L\'evy metric is complete).

\subsection{Proof of Theorem  \ref{th:convmes} (ii) - bounded degree} \label{subsec:th(ii)bounded}

We assume first the following
\begin{itemize}
\item[A".] There exists $\ell \geq 1$ such that for all $n \in \N$
and $i \in [n]$, $\deg(G_n,i) \leq \ell $.
\end{itemize}
With this extra assumption, $\Delta_{n}$ and $\Delta$ are self
adjoint operators and, as in Section \ref{subsec:th(i)bounded}, we
deduce that $\EE U(\oG_n)\Rightarrow \orho$. In particular, we have
\begin{eqnarray*}
\EE[m_{n}(z)]\to \int \la R( z)o,o\ra d\orho[G,o].
\end{eqnarray*}
Hence, in order to prove Theorem \ref{th:convmes} (ii) with the
additional assumption (A"), it is sufficient to prove that, for
all $z \in \C_+$,
\begin{eqnarray}
\label{eq:L2}\lim_{n\to \infty}\EE  \left|m_n(z) -
\EE[m_n(z)]\right|^2 =0.
\end{eqnarray}
Take $z\in \C_+$ with $\Im z> 2\ell$. Notice, that by Vitali's
convergence Theorem, it is sufficient to prove (\ref{eq:L2}) for
all $z$ such that $\Im z> 2 \ell$. By assumption (A"), we get $\left|(\Delta_{n})^k_{ii}
\right|\leq \ell^k$, and we may thus write for any integer $t \geq 1$,
\begin{eqnarray*}
R_n(z)_{ii}(z) &= & -\sum_{k=0}^{t-1}
\frac{(\Delta_n)^k_{ii}}{z^{k+1}}-\sum_{k=t}^\infty
\frac{(\Delta_n)^k_{ii}}{z^{k+1}}\\
&=:& X^{(n)}_i(t) + \epsilon^{(n)}_i(t).
\end{eqnarray*}
Note that $|\epsilon^{(n)}_i(t)|\leq \epsilon(t)
:=\sum_{k = t}^\infty \frac{\ell^k}{(2 \ell)^{k+1}} \leq 2^{-t}$. We define
$\oX^{(n)}_i(t) =
X^{(n)}_i(t)-\EE\left[X^{(n)}_i(t)\right]$. Since
$\left|X^{(n)}_i(t) - R_n(z)_{ii}  \right| =\left|
\epsilon^{(n)}_i(t)  \right| \leq \epsilon(t)$, we have
\begin{eqnarray*}
\left|m_n(z) - \EE[m_n(z)]\right|\leq\left|
  \frac{1}{n}\sum_{i=1}^n \oX^{(n)}_i(t)\right|+2\epsilon(t).
\end{eqnarray*}
Therefore if for all $t \geq 1$, we manage to prove that in $L^2( \PP) $
\begin{equation}
\label{eq:asL} \lim_{n \to \infty} \frac{1}{n} \sum_{i=1} ^ n
\overline  X^{(n)}_{i} (t) = 0,
\end{equation}
then the proof of (\ref{eq:L2}) will be complete. We now prove
(\ref{eq:asL}):
\begin{eqnarray*}
\EE \left( \frac{1}{n} \sum_{i=1} ^ n \overline  X^{(n)}_{i}
(t) \right)^2 & = & \frac{1}{n^2} \EE\left( \sum_{i\neq
j}\oX^{(n)}_{i}
  (t)\oX^{(n)}_{j} (t)+\sum_{i=1}^n \oX^{(n)}_{i} (t)^2\right) \\
  & = & \EE \left( \oX^{(n)}_{o_1} (p )\oX^{(n)}_{o_2} (t) \right),
\end{eqnarray*}
where $(o_1,o_2)$ is a uniform pair of vertices. We then notice
that $\oX^{(n)}_{i} (t)$ is a measurable function of the ball
of radius $t$ and center $i$. Thus, by assumption (R), $\lim_n
\EE \oX^{(n)}_{o_1} ( t )\oX^{(n)}_{o_2} (t)= \lim_n \EE
\oX^{(n)}_{o_1} (t) \EE \oX^{(n)}_{o_2} (t)  = 0$, and
(\ref{eq:asL}) follows. Hence we proved (\ref{eq:L2}) under the
assumption (A").

\subsection{Proof of Theorem  \ref{th:convmes} (ii) - general case} \label{sec:th(ii)}

We now relax assumption (A") by assumption (A). By the same
argument as in Section \ref{sec:th(i)}, we get
$$
\EE L (\mu_{\Delta_{n}}, \mu_{\Delta_{n,\ell}}) \leq \EE  \int ( \deg(G,o) +1) \ind ( \deg(G,o)  > \ell) dU(G_n) [G,o]=: \overline p_{n,\ell}.
$$
By assumptions (R-A), uniformly in $\ell$,
\begin{eqnarray*}
\overline p_{n,\ell}\to \EE \int ( \deg(G,o) +1 ) \ind (
\deg(G,o)  > \ell)d\rho[G,o],
\end{eqnarray*}
where the right-hand term tends to $0$ as $\ell$ tends to infinity. The end of the proof follows by the same argument, since we have
now for $n\geq N(\ell)$, $q\geq 1$,
\begin{eqnarray*}
\EE L (\mu_{\Delta_{n+q}}, \mu_{\Delta_{n}}) \leq  2 \epsilon +\EE  L(\mu_{\Delta_{n+q,\ell}} , \mu_{\Delta_{n,\ell}}).
\end{eqnarray*}

\section{Proof of Theorem \ref{th:adj}} \label{sec:Th2}

\subsection{Proof of Theorem \ref{th:adj} (i)} \label{subsec:unicity}

In this paragraph, we check the existence and the unicity of the
solution of the RDE (\ref{eq:RDE}). Let $\Theta = \N \times
\cH^\infty$, where $\cH^\infty$ is the usual infinite product
space. We define a map $\psi: \Theta \to \cH$ as follows
\begin{eqnarray*}
\begin{array}{cccl}
\psi (n,(h_i)_{i \in \N}) : & \C_+&\to&\C_+\\
&z&\mapsto& -\left(z+\alpha(n+1) + \sum_{i=1}^ {n} h_i (z)
\right)^{-1}.
\end{array}
\end{eqnarray*}
Let $\Psi$ be a map from $\cP(\cH)$ to itself, where $\Psi(P)$ is
the distribution of $\psi(N,(Y_i)_{i \in \N})$, where
\begin{itemize}
\item[(i)] $(Y_i, i\geq 1)$ are independent with distribution $P$;
\item[(ii)] N has distribution $F$; \item[(iii)] the families in
(i) and (ii) are independent.
\end{itemize}
We say $Q \in \cP(\cH)$ is a solution of the RDE (\ref{eq:RDE}) if
$Q= \Psi(Q)$.

\begin{lemma}
\label{le:unicity} There exists a unique measure $Q    \in
\cP(\cH)$ solution of the RDE (\ref{eq:RDE}).
\end{lemma}
\bp Let $\Omega$ be a bounded open set in the half plane $\{ z \in
\C_+ : \Im z \geq \sqrt{\EE N}+1 \}$ (by assumption (A'), $\EE N$
is finite). Let $\cP (\cH)$ be the set of probability measures on
$\cH$. We define the distance on $\cP (\cH)$
$$
W (P,Q) = \inf \EE  \int_\Omega  | X(z) - Y(z) | dz
$$
where the infimum is over all possible coupling of the
distributions $P$ and $Q$ where $X$ has law $P$ and $Y$ has law
$Q$. The fact that $W$ is the distance follows from the fact that
two holomorphic functions equal on a set containing a limit point
are equal. The space $\cP(\cH)$ equipped with the metric $W$ gives
a complete metric space.

Let $X$ with law $P$, $Y$ with law $Q$ coupled so that $W (P,Q) =
\EE  \int_\Omega  | X(z) - Y(z) | dz$. We consider $(X_i,Y_i)_{i
\in \N}$ iid copies of $(X,Y)$ and independent of the variable
$N$. By definition, we have the following
\begin{eqnarray*}
W(\Psi (P), \Psi (Q)) & \leq & \EE  \int_\Omega | \psi(N,(X_i); z) - \psi(N,(Y_i) ; z) | dz \\
 & \leq &  \EE \int_\Omega \left| \frac{ \sum_{i =1}^{N}  X_i (z) - Y_i (z) }{ \left( z + \alpha (N+1) + \sum_{i =1}^{N} X_i (z) \right)  \left( z + \alpha (N+1) + \sum_{i =1}^{N} Y_i (z) \right) }\right| dz \\
 & \leq & \int_\Omega (\Im z)^{-2}  \EE  \sum_{i =1}^{N}  \left| X_i (z) - Y_i (z) \right| dz \\
 & \leq &\EE N  (\inf_{z\in \Omega} \Im z)^{-2}    W(P, Q).
\end{eqnarray*}
Then since $\inf_{z\in \Omega}\Im z >  \sqrt{\EE N} $, $\Psi$ is a
contraction and from Banach fixed point Theorem, there exists a
unique probability measure $Q  $ on $\cH$ such that $\Psi(Q  ) = Q
$. \ep

\subsection{Resolvent of a tree}

In this paragraph, we prove the following proposition.

\begin{proposition}
\label{le:treeD} Let $F_*$ be a distribution with finite mean, and
$(T_n,1)$ be a GWT rooted at $1$ with degree distribution $F_*$
stopped at generation $n$. Let $A(T_n)$ be the adjacency matrix
of $T_n$, and $\Delta (T_n)   = A(T_n)-\alpha D (T_n)$. Let
$R^{(n,T)} (z)= ( \Delta (T_n)  - z I )^{-1}$ and $X^{(n,T)}
(z)= R^{(n,T)}_{11}(z)$. For all $z \in \C_+$, as $n$ goes to
infinity $X^{(n,T)}  (z) $ converges weakly to $X (z)$ defined by
Equation (\ref{eq:X}).
\end{proposition}

We start with the following Lemma which explains where the RDE
(\ref{eq:RDE}) comes from.
\begin{lemma}
\label{le:tree} Let $F$ be a distribution with finite mean, and
$(T_n,1)$  be a GWT rooted at $1$ with offspring distribution $F$
stopped at generation $n$. Let $A(T_n)$ be the adjacency matrix
of $T_n$, and  $\bar{\Delta}(T_n)   =
A(T_n)-\alpha (D(T_n) + V(T_n))$, where $V(T_n)_{11}=1$
and $V(T_n)_{ij}=0$ for all $(i,j)\neq (1,1)$. Let $R^{(n,T)}
(z)= (\bar{\Delta}(T_n)   - z I )^{-1}$ and $Y^{(n,T)} (z)=
R^{(n,T)}_{11}(z)$. For all $z \in \C_+$, as $n$ goes to infinity
$Y^{(n,T)}(z) $ converges weakly to $Y(z)$ given by the RDE
(\ref{eq:RDE}).
\end{lemma}
\bp For simplicity, we omit the superscript $T$ and the variable $z$. 
We order the vertices of $(T_n,1)$ according to a depth-first search in
the tree. We denote by  $N = D^{(n)}_{11}$ the number of
offsprings of the root and by $T^1_n, \cdots, T^{N}_n$ the subtrees of
$T_n\backslash\{1\}$ ordered in the order of the depth-first search.
With this ordering, we obtain a matrix
$\bar{\Delta}(T_n)$ of the following shape:
$$
\begin{pmatrix} 
 0&\begin{matrix}1&0&\cdots\end{matrix}&\begin{matrix}1&0&\cdots\end{matrix}&\cdots&\begin{matrix}1&0&\cdots\end{matrix}\\
\begin{matrix}1\\0\\\vdots\end{matrix}&\bar \Delta ( T^1 _n ) &&& \\
\begin{matrix}1\\0\\\vdots\end{matrix}&&\bar \Delta ( T^2 _n )  && \\
\vdots&&&&\\
\begin{matrix}1\\0\\\vdots\end{matrix}&&&&\bar \Delta ( T^N _n ) 
\end{pmatrix}.$$

We then use a classical of Schur decomposition formula for $\bar{\Delta}(T_n)   - z I $ (see Lemma \ref{le:schur} in Appendix)
\begin{eqnarray}\label{eq:schur1}
Y^{(n)} = -\left(z+\alpha (D(T_n)_{11}+1) + \sum_{2\leq
i,j\leq n} \widetilde{R}^{(n-1)}_{ij}A ( T_n) _{1i}A (T_n) _{1j}
\right)^{-1},
\end{eqnarray}
where $\widetilde{R}^{(n-1)} = (\widetilde{\Delta}^{(n-1)} -
zI )^{-1}$ with $\widetilde{\Delta}^{(n-1)}$ is the matrix
obtained from $\bar{\Delta}(T_n) $ with the first row and column
deleted. 
It follows that $\widetilde{R}^{(n-1)} $ is decomposable in the following diagonal block form
$$
\begin{pmatrix} 
  R ( T^1 _n )  &                          &  &  \\
                             & R  ( T^2 _n ) &     & \\
                              & & \cdots & \\
                              & & & R ( T^N _n ) 
\end{pmatrix},$$
where $R(T^i _n ) =  (\bar \Delta ( T^i _n ) - z I ) ^{-1} $.
In particular, we get
\begin{eqnarray}
\label{eq:schur2}
 \widetilde{R}^{(n-1)}_{ij}A ( T_n) _{1i}A (T_n) _{1j} = 0 \quad \hbox{ if } \quad i \ne j. 
\end{eqnarray}
Indeed, if $A ( T_n) _{1i}A (T_n) _{1j} =1$ then both $i$ and $j$ are
offsprings of $1$ in separate subtrees so that $\widetilde{R}^{(n-1)}_{ij} =0$. 

Since $T_n$ is a Galton Watson tree of depth $n$, the number of offsprings of the root, $N = D^{(n)}_{11}$, has distribution $F$. Moreover the subtrees $T^1_n, \cdots, T^{N}_n$ are iid with common distribution $T_{n-1}$, and are independent of $N$. We now define
$$
Y^{(n-1)}_i  = \left( (\bar \Delta ( T^i _n ) - z I ) ^{-1} \right)_{v_i , v_i} = \widetilde{R}^{(n-1)}_{v_i , v_i}. 
$$
It follows that $(Y^{(n-1)}_1,\cdots, Y^{(n-1)}_N)$ are iid, independent of $N$, with the same common law than $Y^{(n-1)}$. From Equations (\ref{eq:schur1}), (\ref{eq:schur2}), we deduce
\begin{eqnarray*}
Y^{(n)} = -\left(z+\alpha(N+1)+ \sum_{i=1}^N Y^{(n-1)}_i 
\right)^{-1},
\end{eqnarray*}
In other words, with a slight abuse of notation, and identifying a random variable with its distribution, we have $Y^{(n)}
= \Psi(Y^{(n-1)})$, where the mapping $\Psi$ was defined in \S
\ref{subsec:unicity}. The end of the proof follows directly from
Lemma \ref{le:unicity}. \ep

\noindent{\em Proof of Proposition \ref{le:treeD}.} Again, we omit the superscript $T$ and the variable $z$. As above, we use
the decomposition formula:
\begin{eqnarray*}
X^{(n)}   = -\left(z+\alpha D(T_n)_{11} + \sum_{2\leq i,j\leq
n} \widetilde{R}^{(n-1)}_{ij} A ( T_n) _{1i}A (T_n) _{1j}
\right)^{-1},
\end{eqnarray*}
where $\widetilde{R}^{(n-1)} = (\widetilde{\Delta}^{(n-1)} -
zI )^{-1}$ with $\widetilde{\Delta}^{(n-1)}$ is the matrix
obtained from $\Delta (T_n)   $ with the first row and column
deleted. As above, since $T_n$ is a tree $\widetilde{R}^{(n-1)}_{ij}
A ( T_n) _{1i}A (T_n) _{1j}  = 0$ if $i \ne j$, so that we get
\begin{eqnarray*}
X^{(n)}    = -\left(z+\alpha N_*+ \sum_{i=1}^{N_*} Y^{(n-1)}_i
\right)^{-1},
\end{eqnarray*}
where $N_*$ has distribution $F_*$ and $Y^{(n-1)}_i$ are iid
copies of $Y^{(n-1)}$, independent of $N_*$, defined in Lemma
\ref{le:tree}. Proposition \ref{le:treeD} follows easily from
Lemma \ref{le:tree}.\ep

\subsection{Proof of Theorem \ref{th:adj} (ii) - bounded degree}

We first assume that Assumption (A") holds. Let $T$ be a GWT with degree distribution $F_*$ and $T_n$ be the restriction of $T$  to the set of vertices at distance at most $n$ from the root (i.e. $T_n$ is stopped at generation $n$). Let $\Delta = \Delta(T)$ denote the operator associated to $T$, as in Section \ref{sec:graphop}, $\Delta$ is self-adjoint and we may define for all $z \in \C_+$, $R(z) = ( \Delta - zI)^{-1}$.  Then by Theorem
VIII.25(a) in \cite{reesim}, $\Delta(T_n)  $ converges to $\Delta$ in the strong resolvent sense.
Hence, by Proposition \ref{le:treeD} we have $$\la R(z)o,o\ra
\stackrel{d}{=} X   (z).$$

Now, as in \S \ref{sec:graphop}, \ref{subsec:th(ii)bounded},
$\Delta^{(n)}  $ and $\Delta  $ are self adjoint operators and
$\EE U(\oG_n)\Rightarrow \orho$. From Theorem \ref{th:convmes}, we
have
\begin{eqnarray*}
m_{n}(z) \stackrel{L^1}{\longrightarrow} \int \la R(z)o,o\ra
d\orho[G,o].
\end{eqnarray*}
The proof of Theorem \ref{th:adj} (ii) is complete with the extra
assumption (A").

\subsection{Proof of Theorem \ref{th:adj} (ii) - general case}

We now relax assumption (A") by assumption (A). From Theorem
\ref{th:convmes}, it is sufficient to prove that $\lim_n \EE
m_{n} (z)$ converges to $\EE X   (z)$, where $X  $ is defined by
Equation (\ref{eq:X}). By the same argument as in \S
\ref{sec:th(i)}, it is sufficient to prove that, for the weak
convergence on $\cH$,
$$ \lim_{\ell \to \infty} X^{(\ell)}    = X  $$
where $X^{(\ell)}  $ is defined by Equation (\ref{eq:X}) with a
degree distribution $F_{*}^{(\ell)}$ which converges weakly to $F$ as
$\ell$ goes to infinity. This continuity property is established in
Lemma \ref{le:Dn3} (in Appendix). The proof of Theorem \ref{th:adj} is
complete.

\section{Applications and extensions} \label{sec:ext}

\subsection{Weighted graphs}

A weighted graph is a graph $G = (V,E)$ with attached weights on
its edges. As in \S \ref{subsec:conv}, we consider a sequence of
graphs $G_n$ on $[n]$. We define the symmetric matrix $W =
(w_{ij})_{1 \leq i,j \leq n}$, where $(w_{ij})_{1 \leq i \leq j }$
is a sequence of iid real variables, independent of $G_n$ and
$w_{ij} = w_{ji}$. Let $\circ$ denote the Hadamard product (for all $i,j$,
$(A\circ B ) _{ij} = A_{ij} B_{ij}$) and let $T(G_n)$ be the diagonal
matrix whose entry $(i,i)$ is equal to $\sum_{j} w_{ij}
A^{(n)}_{ij}$. We define  $\mu_{\overline \Delta_n}  $ as the spectral
measure of the matrix $\overline \Delta_n    = W \circ A (G_n) - \alpha T(G_n) $. If
$o$ denotes the uniformly picked root of $G_n$, we assume
\begin{itemize}
\item[B.] The sequence of variables $\left(\sum_{k=1}^n
A^{(n)}_{o k }  |w_{o k }|  \right), n \in \N,$
is uniformly integrable.
\end{itemize}
An easy extension of Theorem \ref{th:convmes} is
\begin{theorem}
\begin{itemize}
\item[(i)] Let $G_n = ( [n] , E_n )$ be a sequence of graphs
satisfying assumptions (D-B). Then there exists a probability
measure $\nu  $ on $\R$ such that a.s. $\lim_{n \to \infty}
\mu_{\overline \Delta_n}   = \nu  $. \item[(ii)] Let $G_n = ( [n] , E_n )$ be a
sequence of random graphs satisfying Assumptions (R-B). Then there
exists a probability measure $\nu  $ on $\R$ such that, $\lim_{n
\to \infty} \EE L (\mu_{\overline \Delta_n}   , \nu  ) = 0$.
\end{itemize}
\end{theorem}
The only difference with the proof of  Theorem \ref{th:convmes} appears in \S \ref{sec:th(i)}, for $\ell \in \N$,  the matrix
$\Delta_{n,\ell}$ is now equal to, for $i \neq j$
$$
 (\Delta_{n,\ell})_{ij} = \left\{\begin{array}{ll}

A(G_n)_{ij} w_{ij} & \hbox{ if }\max \{  \sum_{k=1}^n
A^{(n)}_{ik}  |w_{ik}| , \sum_{k=1}^n
A^{(n)}_{jk}  |w_{jk}|\} \leq \ell  \\
0 & \hbox{ otherwise}
\end{array}\right.
$$
and  $(\Delta_{n,\ell})_{ii} = -\alpha \sum_{j \neq
i}  (\Delta_{n,\ell})_{ij}$. The remainder is identical.

We may also state an analog of Theorem \ref{th:adj} for the case
$\alpha = 0$, that is $\overline \Delta_n = W \circ A (G_n)$. We denote by
$s_n$ the Stieltjes transform of $\mu_{\overline \Delta_n}$. Assumption (A') is strengthen into
\begin{itemize}
\item[B'.] Assumption (B) holds, $\sum_k k^2 F_* (k) < \infty$ and $\EE [ w_{12} ^2] < \infty$.
\end{itemize}

The proof of the next
result is a straightforward extension of the proof of Theorem
\ref{th:adj}.

\begin{theorem}
\label{th:adjw} Assume that assumptions (RT-B') hold and $\alpha =
0$ then
\begin{itemize}
\item[(i)] There exists a unique probability measure $P \in
\cP(\cH)$ such that for all $z\in  \mathbb C_+$,
\begin{equation*}
Y (z)   \stackrel{d}{=}  - \left(z  +  \sum_{i=1}^ {N} |w_i |^2
Y_i(z)  \right)^{-1},
\end{equation*}
where $N$ has distribution $F$, $w_i$ are iid copies with
distribution $w_{12}$, $Y$ and $Y_i$ are iid copies with law $P$
and the variables $N,w_i,Y_i$ are independent. \item[(ii)] For all
$z\in \C_+$, $s_n  (z)$ converges as $n$ tends to infinity in
$L^1$ to $ \EE X (z)$, where for all $z \in \C_+$,
$$
X (z)\stackrel{d}{=}  - \left(z  + \sum_{i=1}^ {N_*} |w_i |^2
Y_i(z)   \right)^{-1},
$$
where $N_*$ has distribution $F_*$, $w_i$ are iid copies with
distribution $w_{11}$, $Y_i$ are iid copies with law $P$ and the
variables $N,w_i,Y_i$ are independent.
 \end{itemize}
\end{theorem}

The case $\alpha = 1$ is more complicated: the diagonal term
$T(G_n)$ introduces a dependence within the matrix which breaks
the nice recursive structure of the RDE.

\subsection{Bipartite graphs}

In \S \ref{subsec:tree}, we have considered a sequence of random
graphs converging weakly to a GWT tree. Another important class of
random graphs are the bipartite graphs. A graph $G = (V,E)$ is
bipartite if there exists two disjoint subsets $V^a$, $V^b$, with
$V^a \cup V = V$ such that all edges in $E$ have an adjacent
vertex in $V^a$ and the other in $ V^b$. In particular, if $n$ and
$p$ are the cardinals of $V^a$ and $V^b$,  the adjacency matrix of
a bipartite graph may be written as $\left( \begin{array} {cc} 0 &
M ^* \\ M & 0 \end{array} \right) $ for some $n \times p $ matrix
$M$. The analysis of random bipartite graphs finds strong
motivation in coding theory, see for example Richardson and
Urbanke \cite{richardsonurbanke}. Note also that if $M$ is an $n
\times p $ matrix, the spectrum of $M^*M$ may be obtained from the
spectrum of $\left( \begin{array} {cc} 0 & M ^* \\ M & 0
\end{array} \right) $.   We may thus find another motivation in
sparse statistical problem, see El Karoui \cite{elkaroui}.

The natural limit for random bipartite graphs is the following
Bipartite Galton-Watson Tree (BGWT) with degree distribution
$(F_*,G_*)$ and scale $p \in (0,1)$. The BGWT is obtained from a
Galton-Watson branching process with alternated degree
distribution. With probability $p$, the root has offspring
distribution $F_*$, all odd generation genitors have an offspring
distribution $G$, and all even generation genitors (apart from the
root) have an offspring distribution $F$. With probability $1- p$,
the root has offspring distribution $G_*$, all odd generation
genitors have an offspring distribution $F$, and all even
generation genitors have an offspring distribution $G$.

We now consider a sequence $(G_n)$ of random bipartite graphs
satisfying assumptions $(R-A')$ with weak limit a BGWT with degree
distribution $(F_*,G_*)$ and scale $p \in (0,1)$. The weak
convergence of a natural ensemble of bipartite graphs toward a
BGWT  with degree distribution $(F_*,G_*)$ and scale $p \in (0,1)$
follows from \cite{richardsonurbanke}, $p$ being the proportion of
vertices in $V^a$, $F_*$ the asymptotic degree distribution of
vertices in $V^a$ and $b_*$ the asymptotic degree distribution of
vertices in $V^a$. As usual, we denote by $\mu_{\Delta_n}    $ the
spectral measure of $\Delta_n  $, $m_n  $ is the Stieltjes
transform of $\mu_{\Delta_n}   $. We give
without proof the following theorem which is a generalization of
Theorem \ref{th:adj}.

\begin{theorem}
\label{th:bi} Under the foregoing assumptions,
\begin{itemize}
\item[(i)] There exists a unique pair of probability measures
$(R^a  ,R^b  ) \in \cP(\cH) \times \cP(\cH)$ such that for all
$z\in  \mathbb C_+$,
\begin{eqnarray*}
Y^a (z)    \stackrel{d}{=}   - \left(z  + \alpha (N^a+1) + \sum_{i=1}^ {N^a} Y^b_i(z)  \right)^{-1},  \\
Y^b (z)    \stackrel{d}{=}   - \left(z  + \alpha (N^b+1) +
\sum_{i=1}^ {N^b} Y^a_i(z)  \right)^{-1},
\end{eqnarray*}
where $N^a$ (resp. $N^b$) has distribution $F$ (resp. $G$) and
$Y^a,Y^a_i$ (resp. $Y^b, Y^b_i$) are iid copies with law $R^a  $
(resp. $R^b  $), and the variables  $N^b, Y^a_i,N^a, Y^b_i$ are
independent. \item[(ii)] For all $z\in \C_+$, $m_n   (z)$
converges as $n$ tends to infinity in $L^1$ to $ p \EE X^a   (z) +
(1-p) \EE X^b   (z) $ where for all $z \in \C_+$,
\begin{eqnarray*}
X^a   (z)  \stackrel{d}{=}  - \left(z  + \alpha N^a_*+ \sum_{i=1}^ {N^a_*} Y^b_i(z)   \right)^{-1}, \\
X^b   (z)  \stackrel{d}{=}  - \left(z  + \alpha N^b_*+ \sum_{i=1}^
{N^b_*} Y^a_i(z)   \right)^{-1},
\end{eqnarray*}
where $N^a_*$ (resp. $N^b_*$) has distribution $F_*$ (resp. $G_*$)
and $Y^a_i$ (resp. $Y^b_i$) are iid copies with law $R^a  $ (resp.
$R^b  $), independent of $N^b_*$ (resp. $N^a_*$ ).
 \end{itemize}
\end{theorem}

In the case $\alpha = 0$ and for bi-regular graphs (i.e. BGWT with
degree distribution $(\delta_k,\delta_l)$ and parameter $p$), the
limiting spectral measure is already known and first derived by
Godsil and Mohar \cite{godsilmohar}, see also Mizuno and Sato
\cite{mizunosato} for an alternative proof.

\subsection{Uniform random trees}

The uniformly distributed tree on $[n]$ converges weakly to the
Skeleton tree $T_{\infty}$ which is defined as follows. Consider a
sequence $T_0,T_1, \cdots$ of independent GWT with offspring
distribution the Poisson distribution with intensity $1$ and let
$v_0, v_0, \cdots$ denote their roots. Then add all the edges
$(v_i,v_{i+1})$ for $i \geq 0$. The distribution in $\cG^*$ of the
corresponding infinite tree is the Skeleton tree. See
\cite{aldste} for further properties and Grimmett
\cite{grimmett80} for the original proof  of the  weak convergence
of the uniformly distributed tree on $[n]$ to the Skeleton tree
$T_{\infty}$. 

Let $\mu_{A(G_n)}$ denote the spectral measure of the adjacency
matrix of the random spanning tree $T_n$ on $[n]$ drawn
uniformly (for simplicity of the statement, we restrict ourselves
to the case $\alpha = 0$). We denote by $m_n (z)$ the Stieltjes transform of
$\mu_{A(G_n)}$. 

As an application of Theorems \ref{th:convmes},
\ref{th:adj}, we have the following:

\begin{theorem}Assume $\alpha = 0$.
\begin{itemize}
\item[(i)]
 There exists a unique probability measure $R \in \cP(\cH)$ such that for all $z\in  \mathbb C_+$,
\begin{equation}
\label{eq:RDEut} X (z)   \stackrel{d}{=}  \left(  W(z) ^{-1} - X_1
(z)    \right)^{-1},
\end{equation}
where $W \stackrel{d}{=} ( A(T_0) - z I) ^{-1} _{v_0 v_0} \in \cH$ is the resolvent taken at the root of a GWT with
offspring distribution $Poi(1)$,  $X$ and $X_1$ have law $R$ and
the variables $W$ and $X_1$ are independent. \item[(ii)] For all
$z\in \C_+$, $m_n(z)$ converges as $n$ tends to infinity in
$L^1$ to $ \EE X (z)$.
 \end{itemize}
\end{theorem}

\noindent{\em Sketch of Proof.} The sequence $T_n$ satisfy
(R-A'), thus, from Theorem \ref{th:convmes}, it is sufficient to
show that the resolvent operator $R= (A- z I)^{ -1}$ of the
Skeleton tree taken at the root satisfies the RDE
(\ref{eq:RDEut}). The distribution invariant structure of
$T_\infty$ implies that  $X(z)   \stackrel{d}{=} R(z)_{v_1v_1}$.
The number of offsprings of the root, $v_1$ of $T_1$ is a Poisson
random variable with intensity $1$, say $N$. We denote the
offsprings of $v_1$ in $T_1$ by $v^1_{i_1}, \cdots, v^1_{i_N}$.
From the Schur decomposition formula, we have
$$
R(z)_{v_1v_1} = - \left( z + \widetilde  R (z)_{v_2 v_2} + \sum_{i
= 1}^N R^i (z)_{v^1_i v^1_i} \right)^{-1},
$$
where  $\widetilde  R (z)$ is the resolvent of the infinite tree
obtained by removing $T_1$ and $v_1$ and $R^i (z)$ is the
resolvent of the subtree of the descendants of $v^1_i$ in $T_1$.
Now, by construction $\widetilde  R (z)$ has the same distribution
than $R(z)$, and $R^i (z)$ are independent copies, independent of
$\widetilde R(z)$ of $W(z)$. Thus we obtain
\begin{eqnarray}
X(z)  & \stackrel{d}{=} & - \left( z + X_1(z) + \sum_{i = 1}^N W_i (z) \right)^{-1}, \label{eq:ut1}\\
& \stackrel{d}{=} & - \left(- W(z) ^{-1} +  X_1(z)  \right)^{-1},
\label{eq:ut2}
\end{eqnarray}
where in (\ref{eq:ut2}), we have applied (\ref{eq:RDE}) for GWT
with Poisson offspring distribution. The existence and unicity of
the solution of (\ref{eq:RDEut}) follows from (\ref{eq:ut1}) using
the same proof as in Lemma \ref{le:unicity}. \ep

\section*{Appendix}

For a proof of the next theorem, see e.g. Billingsley \cite{billingsley}. 
\begin{theorem}[Skorokhod Representation Theorem] \label{th:skorokhod}
Let $\mu_n$, $n \in \N$, be a sequence of probability measures on a  complete metric separable space $S$. Suppose that $\mu_n$ converges weakly to a probability measure $\mu$ on $S$. Then there exist random variables $X_n$, $X$ defined on a common probability space $(\Omega, \cF, \PP)$ such that $\mu_n$ is the distribution of $X_n$, $\mu$  is the distribution of $X$, and for every $\omega \in \Omega$, $X_n(\omega)$ converges to $X(\omega)$. 
\end{theorem}

The next lemma is a consequence of Lidskii's inequality. For a proof see Theorem 11.42 in \cite{bai-silverstein-book}
\begin{lemma}[Rank difference inequality]\label{le:rankineq}
\label{le:Dn2} Let $A$, $B$ be two $n \times n$ Hermitian matrices with empirical spectral measures $\mu_A= \frac 1 n \sum_{i=1} ^n \delta_{\lambda_i ( A) }$ and $\mu_B = \frac 1 n \sum_{i=1} ^n \delta_{\lambda_i ( A) }$. Then
$$
L(\mu_A,\mu_B) \leq \frac 1 n \rank( A - B).
$$
\end{lemma}

The next lemma is a standard tool of random matrix theory. 
\begin{lemma}[Schur formula]\label{le:schur}
\label{le:schur} Let $ B=
\begin{pmatrix}
  b_{11} & u^* \\
  u & \tilde B
\end{pmatrix}
$,
denote a $n \times n$ hermitian invertible matrix, $u$ being a vector of dimension $n-1$. Then $$ 
( B^{-1} ) _{11} = \left( b_{11} - u^* \tilde B^{-1} u
\right)^{-1}.
$$
\end{lemma}

\begin{lemma}
\label{le:Dn3} Let $(F_* ^{(n)}), n \in \N,$ be a sequence of
probability measures on $\N$ converging weakly to $F_*$ such that $\sup_n
\sum k F_*^{(n)} (k) < \infty$. Denote by $X^{(n)}  $ and $X  $
the variable defined by (\ref{eq:X}) with degree distribution
$F_*^{(n)}$ and $F_*$ respectively. Then, for the weak convergence
on $\cH$, $\lim_n X^{(n)}   = X  .$
 \end{lemma}
\bp Let $F_*$ and $F'_*$ be two probability measures on $\N$ with
finite mean, and let $d_{TV} (F_*,F'_*) = \sup_{A \subset \N} |
\int_A F_*(dx) - \int_A F'_*(dx) | = 1/2 \sum_k |F_*(k) -
F'_*(k)|$ be the total variation distance. Let $N_*,N'_*,N,N'$
denote variables with law $F_*,F'_*,F,F'$ respectively, and
coupled so that $2 \PP ( N_* \ne N'_*) = d_{TV} (F_*,F'*)$ and $2
\PP ( N \ne N') = d_{TV} (F,F')$ (the existence of these variables
is guaranteed by the coupling inequality). We now reintroduce the
distance defined in the proof of Lemma \ref{le:unicity}. Let
$\Omega$ be a bounded open set in $\C_+$ with an empty
intersection with the ball of center $0$ and radius $\sqrt{\EE
N}+1$. Let $\cP (\cH)$ be the set of probability measures on
$\cH$. We define the distance on $\cP (\cH)$
$$
W (P,P') = \inf \EE  \int_\Omega  | X(z) - X'(z) | dz
$$
where the infimum is over all possible coupling of the
distributions $P$ and $P'$ where $X$ has law $P$ and $X'$ has law
$P'$. With our assumptions, we may introduce the variables $X := X
$ (with law $P$) and $X' := X'  $ (with law $P'$) defined by
(\ref{eq:X}) with degree distribution $F_*$ and $F'_*$
respectively. The proof of the lemma will be complete if we prove
that there exists $C$, not depending on $F_*$ and $F_*'$, such
that
\begin{equation}
\label{eq:lipschitz} W(P,P') \leq C \max ( d_{TV} (F_*,F'_*),
d_{TV} (F,F')) .
\end{equation}
We denote by $Y$ (with law $Q$) and $Y'$ (with law $Q'$) the
variable defined by (\ref{eq:RDE}) with offspring distribution $F$
and $F'$, coupled so that $W (Q,Q') = \EE  \int_\Omega  | Y(z) -
Y'(z) | dz$. We consider $(Y_i,Y'_i)_{i \in \N}$ iid copies of
$(Y,Y')$ and independent of the variable $N_*$. By definition, we
have the following
\begin{eqnarray}
W(P,P') & \leq & \EE \int_\Omega \left| X'(z) - X(z) \right| \ind(N_* \ne N_*') dz \nonumber \\
  & & \quad\quad +  \EE \int_\Omega \left|  \left( z + \alpha N_* + \sum_{i =1}^{N_*} Y_i (z) \right)^{-1}  -  \left( z + \alpha N_* + \sum_{i =1}^{N_*} Y'_i (z) \right)^{-1} \right| dz \nonumber \\
    &\leq & d_{TV} (F_*,F_*') \int_\Omega (\Im z)^{-1} dz +  \int_\Omega (\Im z)^{-2}  \EE  \sum_{i =1}^{N_*}  \left| Y_i (z) - Y'_i (z) \right| dz  \nonumber \\
  & \leq &  d_{TV} (F_*,F_*') \int_\Omega (\Im z)^{-1} dz  + \EE N_*  (\inf_{z\in \Omega} \Im z)^{-2}    W(Q, Q').\label{eq:boundPP}
\end{eqnarray}
We then argue as in the proof of Lemma \ref{le:unicity}. Since
$\Psi(Q) = Q$ and $\Psi'(Q') = Q'$ (where $\Psi'$ is defined as
$\Psi$ with the distribution $F'$ instead of $F$), we get:
\begin{eqnarray*}
W(Q, Q')& = & W(\Psi(Q), \Psi'(Q'))\\
 & \leq & \EE  \int_\Omega | \psi(N,(Y_i); z) - \psi(N',(Y'_i) ; z) | dz \\
  & \leq & d_{TV} (F,F') \int_\Omega (\Im z)^{-1} dz  +  \EE  \int_\Omega | \psi(N,(Y_i); z) - \psi(N,(Y'_i) ; z) | dz \\
  & \leq & d_{TV} (F,F') \int_\Omega (\Im z)^{-1} dz  + \EE N  (\inf_{z\in \Omega} \Im z)^{-2}    W(Q, Q').
\end{eqnarray*}
Then since $\inf_{z\in \Omega}\Im z >  \sqrt{\EE N} $, we deduce
that
$$  W (Q, Q') \leq d_{TV} (F,F')   \frac{\int_\Omega (\Im z)^{-1} dz}{1 - \EE N  (\inf_{z\in \Omega} \Im z)^{-2}  }.$$
This last inequality, together with (\ref{eq:boundPP}), implies
(\ref{eq:lipschitz}). \ep

\section*{Acknowledgment}
The authors thank Noureddine El Karoui for fruitful discussions
and his interest in this work.

\bibliographystyle{abbrv}
\bibliography{../mat}

\end{document}